\documentclass[12pt]{amsart}
\usepackage{amssymb,amsmath,amsthm}

\numberwithin{equation}{section}
\begin{document}

\theoremstyle{plain}
\newtheorem{theorem}{Theorem}[section]
\newtheorem{lemma}[theorem]{Lemma}
\newtheorem{proposition}[theorem]{Proposition}
\newtheorem{corollary}[theorem]{Corollary}
\newtheorem{conjecture}[theorem]{Conjecture}
\newtheorem*{main}{Main Theorem}

\def\mod#1{{\ifmmode\text{\rm\ (mod~$#1$)}
\else\discretionary{}{}{\hbox{ }}\rm(mod~$#1$)\fi}}

\theoremstyle{definition}
\newtheorem*{definition}{Definition}

\theoremstyle{remark}
\newtheorem*{example}{Example}
\newtheorem*{remark}{Remark}
\newtheorem*{remarks}{Remarks}

\newcommand{\qq}{{\mathbb Q}}
\newcommand{\rr}{{\mathbb R}}
\newcommand{\nn}{{\mathbb N}}
\newcommand{\zz}{{\mathbb Z}}
\newcommand{\ch}{\raisebox{2pt}{$\chi$}}
\newcommand{\al}{\alpha}
\newcommand{\be}{\beta}
\newcommand{\ga}{\gamma}

\newcommand{\ep}{\epsilon}
\newcommand{\la}{\lambda}
\newcommand{\de}{\delta}
\newcommand{\Del}{\Delta}
\title{Regularity properties of the Stern enumeration of the rationals}

\author{Bruce Reznick}
\address{Department of Mathematics, University of 
Illinois at Urbana-Champaign, Urbana, IL 61801} 
\email{reznick@math.uiuc.edu}

\begin{abstract} The Stern sequence $s(n)$ is defined by $s(0) = 0,
s(1) = 1$, $s(2n) = s(n)$, $s(2n+1) = s(n) + s(n+1)$. Stern showed
in 1858 that $gcd(s(n),s(n+1)) = 1$, and that for every pair of relatively
prime positive integers $(a,b)$ there exists a unique $n\ge 1$ with
$s(n) = a$ and $s(n+1) = b$. We show that  in a strong sense, the
average value of $\frac{s(n)}{s(n+1)}$ is $\frac 32$, and that
for $d \ge 2$, $(s(n), s(n+1))$
is uniformly distributed among all feasible pairs of congruence
classes modulo $d$. More precise results are presented for $d=2$ and 3.
\end{abstract}
\date{\today}

\maketitle

\section{Introduction and History}
In 1858, M. A. Stern \cite{ste}
defined the  {\it diatomic  array}, an 
unjustly neglected mathematical construction. It is a
Pascal triangle with memory: each row is created by
inserting the sums of pairs of consecutive elements into the previous row.
\begin{equation}
\begin{aligned}
&a\quad b \\
&a\quad a+b\quad b \\
&a\quad 2a+b\quad a+b\quad a+2b\quad b \\
&a\quad 3a+b\quad 2a+b\quad 3a+2b\quad a+b\quad 2a+3b\quad a+2b\quad
a+3b\quad b \\ 
&\qquad \vdots \\
\end{aligned}
\end{equation}

When $(a,b) = (0,1)$, it is easy to see that each
row of the diatomic array repeats as the first half of the next row
down. The resulting infinite  {\it Stern sequence} can also be
defined recursively by:
\begin{equation}
s(0) = 0,\ s(1) = 1,\qquad s(2n) = s(n),\  s(2n+1) = s(n) + s(n+1).
\end{equation}
Taking $(a,b) = (1,1)$ in (1.1), we obtain blocks of $(s(n))$ for
$2^r \le n \le 2^{r+1}$. Although $s(2^r)=1$ is repeated at the ends,
each pair $(s(n),s(n+1))$ appears below exactly once as a consecutive pair
in a row: 
\begin{equation}
\begin{aligned}
& (r=0) &\qquad  &1\quad 1 \\
& (r=1) &\qquad   &1\quad 2\quad 1 \\
& (r=2) &\qquad   &1\quad 3\quad 2\quad 3\quad 1 \\
& (r=3) &\qquad   &1\quad 4\quad 3\quad 5\quad 2\quad 5\quad 3\quad
  4\quad 1 \\  
&\qquad \vdots \\
\end{aligned}
\end{equation}
Mirror symmetry (or an easy induction) implies that
for $0 \le k \le 2^r$, we have
\begin{equation}
s(2^r + k) = s(2^{r+1} - k).
\end{equation}

In his original paper, Stern proved that for all $n$,
\begin{equation}
\gcd(s(n),s(n+1)) = 1;
\end{equation}
moreover, for every pair of positive relatively prime integers
$(a,b)$, there is a unique $n$ so that $s(n) = a$ and $s(n+1) =
b$. Stern's discovery predates Cantor's proof of the countability of
$\qq$ by fifteen years.
 This property of the Stern sequence has been recently made
explicit and discussed in \cite{cw}. Another enumeration of the
positive rationals involves the {\it Stern-Brocot array}, which also
predates Cantor; see \cite{gkp}, pp. 116--123, 305--306. This was used
by Minkowski in defining his  $?$-function; see \cite{mink}. The
Stern sequence and Stern-Brocot array make brief appearances in
Dickson's {\it History}, see \cite{dic}, pp. 156, 426.
Apparently, de Rham \cite{dr} was the first to consider the
sequence $(s(n))$ {\it per se}, attributing the term ``Stern
sequence'' to Bachmann
\cite{bach}, p. 143, who had only considered the array.
The Stern sequence has recently arisen as well in the discussion of
2-regular sequences \cite{as} and the Tower of Hanoi graph
\cite{hkmpp}. Some other Stern identities and a large bibliography
relating to the Stern sequence are given in \cite{urb}. 
A further discussion of the Stern sequence will be found in \cite{rez2}. 

Let 
\begin{equation}
t(n) = \frac {s(n)}{s(n+1)}.
\end{equation}
Here are blocks of $(t(n))$, for $2^r \le n <2^{r+1}$ for small $r$: 
\begin{equation}
\begin{aligned}
& (r=0) &\qquad &  \tfrac 11 \\  
& (r=1) &\qquad & \tfrac 12 \quad \tfrac 21 \\
& (r=2) &\qquad & \tfrac 13 \quad  \tfrac 32 \quad \tfrac 23 \quad
  \tfrac 31 \\ 
& (r=3) &\qquad & \tfrac 14 \quad \tfrac 43 \quad \tfrac 35 \quad \tfrac 52
\quad \tfrac 25 \quad  \tfrac 53 \quad \tfrac 34 \quad \tfrac 41 \\
&\qquad \vdots \\
\end{aligned}
\end{equation}

In Section 3, we shall show that 
\begin{equation}
\sum_{n=0}^{N-1} t(n) =   \frac {3N}2 + {\mathcal O}(\log^2 N),
\end{equation}
so the ``average'' element in the Stern enumeration of $\mathbb Q_+$ is
$\frac 32$.

For a fixed integer $d \ge 2$, let
\begin{equation}
S_d(n) := (s(n) \text{ mod }d, s(n+1)\text{ mod }  d)
\end{equation}
and let
\begin{equation}
{\mathcal S}_d = \{(i \text{ mod }d, j\text{ mod }d): \gcd(i,j,d) = 1\}.
\end{equation}
It follows from (1.5) that $S_d(n) \in {\mathcal S}_d$ for all $n$. 
In Section 4, we shall show that for each $d$,  the sequence
$(S_d(n))$ is uniformly distributed on
${\mathcal S}_d $, so the ``probability'' that
$s(n) \equiv i \mod d$ can be explicitly computed. 
More precisely, let
\begin{equation}
T(N;d,i) = \left\vert\{n: 0 \le n < N\ \& \  s(n) \equiv i \text{ mod }d
\}\right\vert.
\end{equation}
Then there exists $\tau_d < 1$ so that
\begin{equation}
T(N;d,i) = r_{d,i} N +  {\mathcal O}(N^{\tau_d}),
\end{equation}
where
\begin{equation}
r_{d,i} =  \frac 1d
\cdot \prod_{p | i, p | d} \frac {p}{p+1}  
\cdot \prod_{p \nmid i, p | d} \frac {p^2}{p^2-1}.
\end{equation}
In particular, the probability that $s(n)$ is a multiple of $d$ is
$I(d)^{-1}$, where
\begin{equation}
I(d) = d \prod_{p\ | \ d} \frac{p+1}p \in \mathbb N.
\end{equation}

In Section 5, we present more specific information for the cases $d =
2$ and 3. It is an easy induction
 that $s(n)$ is even if and only if $n$ is a
multiple of 3, so that $\tau_2 = 0$. We show that $\tau_3 =
\frac 12$ and give an explicit formula for $T(2^r;3,0)$, as well as a
recursive description of those $n$ for which $3 \ | \ s(n)$. We also prove
that, for all $N\ge 1$, $T(N;3,1) - T(N;3,2) \in \{0,1,2,3\}$. 

It will be proved in  \cite{rez2} that
\begin{equation}
 T(2^r;4,0) = T(2^r;5,0), \qquad 
T(2^r;6,0) = T(2^r;9,0) = T(2^r;11,0);
\end{equation}
we conjecture that $T(2^r;22,0) = T(2^r;27,0)$. (The latter is
true for $r \le 19$.) These 
exhaust the possibilities for $T(2^r;N_1,0) = T(2^r;N_2,0)$ with
$N_i \le 128$.
Note that $I(4) = I(5) = 6$, $I(6)=I(8) =I(9) = I(11) = 12$ and $I(22)=I(27)
= 36$. However, $T(2^r;8,0)\neq T(2^r;6,0)$, so there is more 
than just asymptotics at work.

\section{Basic facts about the Stern sequence}
We formalize the definition of the diatomic array. Define
$Z(r,k) = Z(r,k;a,b)$ recursively for $r 
\ge 0$ and $0 \le k \le 2^r$ by: 
\begin{equation}
\begin{gathered}
Z(0,0) = a, \quad  Z(0,1) = b; \\ Z(r+1,2k) = Z(r,k), \quad
Z(r+1,2k+1) = Z(r,k) + Z(r,k+1).
\end{gathered}
\end{equation}
The following lemma follows from (1.2), (2.1) and a simple induction.
\begin{lemma}
For $0 \le k \le 2^r$, we have
\begin{equation}
Z(r,k;0,1) = s(k).
\end{equation}
\end{lemma}
Lemma 2.1 leads directly to a general formula for  the diatomic array.
\begin{theorem}
For $0 \le k \le 2^r$, we have
\begin{equation}
Z(r,k;a,b) = s(2^r-k)a + s(k)b. 
\end{equation}
\end{theorem}
\begin{proof}
Clearly, $Z(r,k;a,b)$ is linear in $(a,b)$ and it also satisfies a mirror
symmetry
\begin{equation}
Z(r,k;a,b) = Z(r,2^r-k;b,a)
\end{equation}
for $0 \le k \le 2^r$, c.f. (1.4). Thus,
\begin{equation}
Z(r,k;a,b) = a Z(r,k;1,0) + bZ(r,k;0,1) = aZ(r,2^r-k;0,1) + bZ(r,k;0,1).
\end{equation}
The result then follows from Lemma 2.1.
\end{proof}
The diatomic array contains a self-similarity: any two consecutive entries
in any row determine the corresponding  portion of the succeeding rows. More
precisely, we have a relation whose simple inductive proof is omitted,
and which immediately leads to the iterated generalization of (1.2).
\begin{lemma}
If $0 \le k \le 2^r$ and $0 \le k_0 \le 2^{r_0}-1$, then
\begin{equation}
Z(r+r_0,2^{r}k_0 + k;a,b) =
Z(r,k;Z(r_0,k_0;a,b),Z(r_0,k_0+1;a,b)).
\end{equation}
\end{lemma}
\begin{corollary}
If $n \ge 0$ and $0 \le k \le 2^r$, then
\begin{equation}
s(2^rn + k) = s(2^r-k)s(n) + s(k)s(n+1).
\end{equation}
\end{corollary}
\begin{proof}
Take $(a,b,k_0,r_0) = (0,1,n, \lceil \log_2 (n+1) \rceil)$ in  Lemma
2.3, so that $k_0 < 2^{r_0}$, and then apply Theorem 2.2.
\end{proof}

We turn now to $t(n)$. Clearly, $t(2n) <
1 \le t(2n+1)$ for all $n$; after a little algebra, (1.2) implies
\begin{equation}
t(2n) =  \cfrac 1{1 + \cfrac 1{t(n)}}\ , 
\qquad t(2n+1) = 1 + t(n).
\end{equation}
The mirror symmetry (1.4) yields two other formulas which are
evident in (1.7):
\begin{equation}
t(2^r + k)t(2^{r+1}-k-1) = 1,
\end{equation}
for $0 \le k \le 2^r-1$, which follows from 
\begin{equation}
t(2^{r+1}-k-1) = \frac
{s(2^{r+1}-k-1)}{s(2^{r+1}-k)} = \frac {s(2^r + k+1)}{s(2^r+k)} =
\frac 1{t(2^r+k)};
\end{equation}
 and
\begin{equation}
t(2^r + 2\ell)+ t(2^{r+1}-2\ell-2) = 1,
\end{equation}
for $r \ge 1$ and $0 \le 2\ell \le 2^r-2$, which follows from 
\begin{equation} 
\frac{s(2^r+2\ell)}{s(2^r+2\ell+1)} +
\frac{s(2^{r+1}-2\ell-2)}{s(2^{r+1}-2\ell-1)}
= \frac{s(2^r+2\ell)}{s(2^r+2\ell+1)} +
\frac{s(2^r+2\ell+2)}{s(2^r+2\ell+1)},
\end{equation}
since $s(2m) + s(2m+2) = s(2m+1)$. 

Although we will not use it directly here, we mention a simple
closed formula for $t(n)$, and hence for $s(n)$. Stern had already
proved that if $2^r \le n < 2^{r+1}$, then the sum of the
denominators in the continued fraction representation of $t(n)$ is
$r+1$; this is clear from (2.8). Lehmer \cite{leh} gave an exact
formulation, of which the following is a variation.
 Suppose $n$ is odd and $[n]_2$, the binary
representation of $n$,  consists of a block of $a_1$ 1's, followed by
$a_2$ 0's, $a_3$ 1's, etc, 
ending with $a_{2v}$ 0's and $a_{2v+1}$ 1's, with $a_j \ge 1$. (That
is, $n = 2^{a_1 + \cdots + a_{2v+1}} - 2^{a_2 + \cdots + a_{2v+1}} \pm
\cdots \pm 
2^{a_{2v+1}} - 1$.) Then
\begin{equation}
t(n) = \frac{s(n)}{s(n+1)} = \frac pq = a_{2v+1} + \cfrac 1{a_{2v} +
  \cfrac 1{\dots + \cfrac 1{a_1}}}\ .
\end{equation}
Conversely, if $\frac pq > 1$ and (2.13) gives its presentation as a
simple continued fraction with an odd number of denominators, then the
unique $n$ with $t(n) = \frac pq$ has the binary representation
described above.  (If $n$ is even or $\frac pq < 1$,  apply (2.9) first.) 

The {\it Stern-Brocot array} is named after the
clockmaker Achille Brocot, who used it  \cite{broc} in 1861
 as the basis of a gear table; see also \cite{hayes}. This array
 caught the attention of several French number theorists, and is
 discussed in \cite{luc}. It is formed by applying the diatomic rule
 to numerators and denominators simultaneously:
\begin{equation}
\begin{aligned}
& (r=0) &\qquad &  \tfrac 01 \quad \tfrac 10\\  
& (r=1) &\qquad & \tfrac 01 \quad \tfrac 11 \quad \tfrac 10 \\
& (r=2) &\qquad & \tfrac 01 \quad  \tfrac 12 \quad \tfrac 11 \quad \tfrac
  21\quad \tfrac 10  \\ 
& (r=3) &\qquad & \tfrac 01 \quad \tfrac 13 \quad \tfrac 12 \quad \tfrac 23
\quad \tfrac 11 \quad  \tfrac 32 \quad \tfrac 21 \quad \tfrac 31\quad
\tfrac 10\\
&\qquad \vdots \\  
\end{aligned}
\end{equation}

This array is not quite the same as (1.7).
If $\frac ab$ and $\frac cd$ are consecutive in the $r$-th
row, then they repeat in the $(r+1)$-st row, separated by $\frac
{a+c}{b+d}$. It is easy to see that the elements of the $r$-th row are
$\frac{s(k)}{s(2^r-k)}$, $0 \le k 
\le 2^r$. It is also easy to show that the elements of each row are
increasing, and moreover, that they share a property with the Farey
sequence.
\begin{lemma} For $0 \le k \le 2^r - 2$,
\begin{equation}
\frac{s(k+1)}{s(2^r-k-1)} - \frac{s(k)}{s(2^r-k)} =\frac
1{s(2^r-k)s(2^r-k-1)}. 
\end{equation}
That is,
\begin{equation}
s(k+1)s(2^r-k)-s(k)s(2^r - k-1) = 1.
\end{equation}
\end{lemma}
This lemma has a simple proof by induction, which can be found in
\cite{luc}, p.467  and \cite{gkp}, p.117.    

The ``new'' entries 
in the $(r+1)$-st row of (2.14) are a permutation of the $r$-th row of
(1.7). The easiest way to express the connection (see \cite{rez2}) 
for rationals $\frac pq > 1$ is that if $0 < k < 2^r$ is odd, then
\begin{equation}
\frac pq = \frac{s(2^r+k)}{s(2^r-k)} = \frac {s(\overleftarrow{2^r+k})}
{s(\overleftarrow{2^r+k}+1)},
\end{equation}
where $\overleftarrow n$ denotes the integer so that $[n]_2$ and
$[\overleftarrow n]_2$ are the reverse of each other. If $\frac pq <
1$, then apply mirror symmetry to the instance of (2.17) which holds
for $\frac qp$. 

The  Minkowski $?$-function can be defined using the first half
of the rows of (2.14). For odd $\ell$, $0 \le \ell \le 2^r$, 
\begin{equation}
?\left( \frac{s(\ell)}{s(2^{r+1}-\ell)} \right) = \frac \ell{2^r}.
\end{equation}
This gives a strictly increasing map from $\mathbb Q \cap [0,1]$ to
the dyadic rationals in $[0,1]$, which extends to a continuous
strictly increasing map from $[0,1]$ to itself, taking quadratic
irrationals to non-dyadic rationals.

Finally, suppose $N$ is a positive integer, written as 
\begin{equation}
N = 2^{r_1} + 2^{r_2} + \cdots + 2^{r_v}, \qquad r_1 > r_2 > \dots > r_v.
\end{equation}
We shall define 
\begin{equation}
N_0 = 0;\quad  N_j = 2^{r_1} +\cdots + 2^{r_j}\text{ for } j = 
1, \dots, v. 
\end{equation}
Further,  for $1 \le j \le v$, let $M_j = 2^{-r_j}N_{j+1}$, so that
\begin{equation}
N_j = N_{j-1} + 2^{r_j} = 2^{r_j}(M_j+1) = 2^{r_{j-1}}M_{j-1}.
\end{equation} 
and, for $a < b \in \mathbb Z$, let 
\begin{equation}
[a,b):= \{k \in \mathbb Z : a \le k < b \}.
\end{equation}
Our proofs will rely  on the observation that
\begin{equation}
[0,N) = \bigcup_{j=0}^{v-1} [N_j,N_{j+1}) = \bigcup_{j=1}^{v}
    [2^{r_j}M_j,2^{r_j}(M_j+1)),
\end{equation}
where the above unions are disjoint, so that, formally, 
\begin{equation}
\sum_{n=0}^{N-1}
= \sum_{j=0}^{v-1} \sum_{n=N_j}^{N_{j+1}-1} =
 \sum_{j=1}^{v} \sum_{n=2^{r_j}M_j}^{2^{r_j}(M_j+1)-1}.
\end{equation}

\section{The Stern-Average Rational}
We begin by looking at the sum of $t(n)$ along the rows of (1.7). Let
\begin{equation}
A(r) = \sum_{n=2^r}^{2^{r+1}-1} t(n)\qquad{\text{and}}\qquad
\tilde A(r) = \sum_{n=0}^{2^{r}-1} t(n) = \sum_{i=0}^{r-1} A(i). 
\end{equation}
\begin{lemma}
For $r \ge 0$,
\begin{equation}
A(r) = \frac 32 \cdot 2^r - \frac 12\qquad{\text{and}}\qquad \tilde
A(r) = \frac 32 \cdot 2^r -\frac {r+3}2. 
\end{equation}
\end{lemma}
\begin{proof}
First note that $A(0) = t(1) = \frac 11 = \frac 32 - \frac 12$. Now
observe that for $r \ge 0$, 
\begin{equation}
\begin{gathered}
A(r+1) = \sum_{j=0}^{2^{r+1}-1} t(2^{r+1} + j) =  \sum_{k=0}^{2^r-1}
t(2^{r+1} + 2k) + \sum_{k=0}^{2^r-1} t(2^{r+1} + 2k+1).
\end{gathered}
\end{equation}
Using (2.11) and (2.8), we can simplify this summation:
\begin{equation}
 \sum_{k=0}^{2^r-1} t(2^{r+1} + 2k) = \frac 12 \left( 
 \sum_{k=0}^{2^r-1} t(2^{r+1} + 2k) + t(2^{r+2}-2k-2) \right) 
 = 2^{r-1},
\end{equation}
and
\begin{equation}
 \sum_{k=0}^{2^r-1} t(2^{r+1} + 2k+1) =  \sum_{k=0}^{2^r-1}\bigl(1 +
 t(2^r+k)\bigr) = 2^r + A(r).
\end{equation}
Thus, $A(r+1) = 2^{r-1}+2^r+A(r)$, and the formula for $A(r)$ is
established by induction. This also immediately implies the formula for
$\tilde A(r)$. 
\end{proof}
\begin{lemma}
If $m$ is even, then
\begin{equation}
\tilde A(r) \le \sum_{k=0}^{2^r-1} t(2^rm + k) < A(r).
\end{equation}
 \end{lemma}
\begin{proof}
For fixed $(k,r)$, let
\begin{equation}
\Phi_{k,r}(x) = \frac {s(2^r-k) x + s(k)}{s(2^r-(k+1))x + s(k+1)}.
\end{equation}
Then it follows from (2.16) that 
\begin{equation}
\Phi_{k,r}'(x)  = \frac 
{s(k+1)s(2^r-k)-s(k)s(2^r - k-1)}{(s(2^r-(k+1))x + s(k+1))^2} > 0.
\end{equation}
Using (2.7), we see that
\begin{equation}
\begin{gathered}
t(2^rm + k) = \frac{s(2^rm + k)}{s(2^rm + k+1)} =
\frac{s(2^r-k)s(m) + s(k)s(m+1)}{s(2^r-k-1)s(m) + s(k+1)s(m+1)}
\\=\Phi_{k,r}\left(\frac {s(m)}{s(m+1)} \right) = \Phi_{k,r}(t(m)).
\end{gathered}
\end{equation}
Since $m$ is even, $0 \le t(m) < 1$; monotonicity then implies that
\begin{equation}
t(k) = \Phi_{k,r}(0) \le t(2^rm + k) < \Phi_{r,k}(1) = t(2^r+k).
\end{equation}
 Summing (3.10) on $k$ from 0 to $2^r-1$ gives (3.6).
\end{proof}
We use these estimates to establish (1.8).
\begin{theorem}
If $2^r \le N < 2^{r+1}$, then
\begin{equation}
\frac{3N}2 -  \frac{r^2+7r+6}4 \le  \sum_{n=0}^{N-1} t(n) <  \frac{3N}2 - 
\frac 12.
\end{equation}
\end{theorem}
\begin{proof}
Recalling (2.24), we apply Lemma 3.2 for each $j$, with $r = r_j$ and
$m = M_j$, so that 
\begin{equation}
 \frac 32 \cdot 2^{r_j} - \frac {r_j+3}2 \le 
 \sum_{n=N_{j-1}}^{N_j-1}  t(n)
 <  \frac 32\cdot  2^{r_j} - \frac 12.
\end{equation}
 After summing on $j$, we find that
\begin{equation}
 \frac {3N}2 - \frac {r_1 + \dots + r_v + 3v}2 \le 
\sum_{n=0}^{N-1} t(n) <  \frac {3N}2 - \frac {v}2.
\end{equation}
To obtain (3.11), note that
$\sum r_j + 3v \le \frac{r(r+1)}2 + 3r+3 =
\frac {r^2+7r+6}2$. 
\end{proof}
\begin{corollary}
\begin{equation}
 \sum_{n=0}^{N-1} t(n) =  
\frac {3N}2 +  \mathcal O\left( \log^2N \right).
\end{equation}
\end{corollary}
Since $t(2^r-1) = \frac r1$, the true error term is at least $\mathcal
O(\log N)$. Numerical computations using Mathematica suggest
 that $\log^2N$ can be replaced by $\log N \log\log N$. 
It also seems that, at least for small fixed
positive integers $t$, 
\begin{equation}
\al_t:= \lim_{N \to \infty} \frac 1N \sum_{n=0}^{N-1} \frac {s(n)}{s(n+t)}
\end{equation}
exists. We have seen that $\al_1 = \frac 32$, and if they exist, the
evidence suggests that
 $\al_2\approx 1.262$, $\al_3 \approx 1.643$ and $\al_4 \approx 1.161$. We
are unable to present an explanation for these specific numerical values.

\section{Stern Pairs, mod $d$}

We fix $d \ge 2$ with prime factorization
$d = \prod p_\ell^{e_\ell}$, $e_\ell \ge 1$, and recall the definitions
of $\mathcal S_d$ and $S_d(n)$ from (1.9) and (1.10). Let 
\begin{equation}
N_d = \left\vert \mathcal S_d \right\vert, 
\end{equation}
and for $0 \le i < d$, let
\begin{equation}
N_d(i)= \left\vert \{j\text{ mod } d: (i \text{ mod } d, j \text{ mod
} d) \in 
\mathcal S_d\} \right\vert.
\end{equation}

We now give two lemmas whose proofs rely on the Chinese Remainder Theorem.
\begin{lemma}
The map $S_d: \nn \to \mathcal S_d$ is surjective. 
\end{lemma}
\begin{proof}
Suppose $\al =(i,j) \in
\mathcal S_d$ with $0 \le i, j \le d-1$.
 We shall show that there exists $w \in \mathbb N$ so that
$\gcd(i, j+wd) = 1$. Consequently, there exists $n$ with  $s(n) = i$ and
$s(n+1) = j +  wd$, so that $S_d(n) = \al$.

Write $i= \prod_\ell {q_\ell^{f_\ell}}$, $f_\ell \ge 1$, with $q_\ell$ prime.
 If $q_\ell\ |\ j$, then    $q_\ell\ \nmid \
d$. There exists  $w\ge 0$ so that $w \equiv d^{-1} \mod
{q_\ell^{f_\ell}}$ if 
$q_\ell\ |\ j$ and  $w \equiv 0 \mod {q_\ell^{f_\ell}}$ if
$q_\ell\ \nmid\ j$. Then  $j + wd \not\equiv 0 \mod
{q_\ell^{f_\ell}}$ for all $\ell$, so no prime dividing $i$
divides $j + wd$, as desired.
\end{proof}
\begin{lemma} 
For $0 \le i \le d-1$,
\begin{equation}
N_d = d^2 \prod_{\ell} \frac {p_\ell^2-1}{p_\ell^2} \qquad
\text{and} \qquad
N_d(i) = d \prod_{p_\ell\ |\ i} \frac {p_\ell-1}{p_\ell}. 
\end{equation}
\end{lemma}
\begin{proof}
To compute $N_d$, we use the Chinese Remainder Theorem by
counting the choices for $(i \text{ mod } p_\ell^{e_\ell}, j \text{ mod
}  p_\ell^{e_\ell})$ for each $\ell$. Missing are those 
 $(i,j)$ in which $p_\ell$ divides both $i$ and $j$, and so
the total number of classes is $(p_\ell^{e_\ell} -
p_\ell^{e_\ell-1})^2$ for each $\ell$.

Now fix $i$.
If $p_\ell \ |\ i$, then $(i,j) \in \mathcal S_d$ if and only
if $p_\ell \nmid j$; if $p_\ell \nmid i$, then there is no restriction on
$j$. Thus,  there are either $p_\ell^{e_\ell} -
p_\ell^{e_\ell-1}$ or $p_\ell^{e_\ell}$ choices for $j$, respectively.
\end{proof}

Suppose $\al = (i,j) \in \mathcal
S_d$; let $L(\al):= (i, i+j)$ and $R(\al) = (i+j,j)$, where $i+j$ is
reduced mod $d$ if necessary. Then $L(\al), R(\al) \in \mathcal S_d$ and
the following lemma is immediate.
\begin{lemma}
For all $n$, we have  $S_d(2n) = L(S_d(n))$ and $S_d(2n+1) =
R(S_d(n))$.
\end{lemma}

We now define the directed graph $\mathcal G_d$ as follows. 
The vertices of $\mathcal G_d$ are the elements of $\mathcal S_d$.  
The edges of $\mathcal G_d$ consist of $(\al, L(\al))$ and $(\al,
R(\al))$ where $\al \in \mathcal S_d$. 
Iterating, we see
that $L^k(\al) = (i, i+kj)$ and $R^k(\al) = (i + kj, j)$, so that
$L^d = R^d = id$, and $L^{-1} = L^{d-1}$ and $R^{-1} = R^{d-1}$.
Thus, if $(\al, \be)$ is an edge of $\mathcal G_d$, then
there is a walk of length $d-1$ from $\be$ to $\al$.

Each vertex of $\mathcal G_d$ has
out-degree two; since $(L^{-1}(\al), \al)$ and $(R^{-1}(\al), \al)$
are edges, each vertex has in-degree two as well.
Let $M_d = [m_{\al(d)\be(d)}] = [m_{\al\be}]$ denote the  adjacency
 matrix for $\mathcal G_d$: $M_d$ is
 the $N_d \times N_d$ 0-1 matrix so that  
$m_{\al L(\al)} = m_{\al R(\al)} = 1$, with other entries equal to 0.
For a positive integer $r$, write 
\begin{equation}
M_d^r = [m_{\al\be}^{(r)}];
\end{equation}
then $m_{\al\be}^{(r)}$ is the number of walks of length $r$
from $\al$ to $\be$.
Finally, for $\ga \in \mathcal S_d$, and integers $U_1 < U_2$, let
\begin{equation}
B(\ga;U_1,U_2) = \left\vert \{m: U_1 \le m < U_2\ \&\ S_d(m) = \ga \}
\right\vert 
\end{equation}
The following is essentially  equivalent to Lemma 2.3.
\begin{lemma}
Suppose $\al = S_d(m)$, $\be \in \mathcal S_d$ and $r \ge 1$. Then 
$B(\be;2^rm,2^r(m+1)) =m_{\al\be}^{(r)} $ is equal to the number of
  walks of length $r$ in $\mathcal G_d$ from $\al$ to $\be$. 
\end{lemma}
\begin{proof}
The walks of length 1 starting from $\al$ are
 $(\al, L(\al))$ and $(\al, R(\al))$; that is,
$(S_d(n),S_d(2n))$ and $(S_d(n),S_d(2n+1))$. The rest is an easy induction.
\end{proof}
\begin{lemma}
For sufficiently large $N$, $m_{\al\be}^{(N)} > 0$ for all $\al, \be$.
\end{lemma}
\begin{proof}
Let $\al_0 = (0,1) = S_d(0)$. Note that $L(\al_0) = \al_0$, hence if
there is a walk of length $w$ from $\al_0$ to $\ga$, then there are
such walks of every length $\ge w$. 
By Lemma 4.1, for each $\al \in \mathcal S_d$, there
exists $n_\al$ so that $S_d(n_\al) = \al$. Choose $r$ sufficiently
large that  $n_\al < 2^r$ for all $\al$. Then by Lemma 4.4, for every
$\ga$, there is a walk of length $r$ from $\al_0$ to $\ga$, and so 
there is a  walk of length $(d-1)r$ from  $\ga$ to $\al_0$. Thus, for
any $\al, \be \in \mathcal S_d$, there is at least one walk of length
$dr$ from 
$\al$ to $\be$ via $\al_0$. 
\end{proof}

We need a version of Perron-Frobenius. Observe that $A_d = \frac 12
M_d$  is doubly
stochastic and the entries of $A_d^N = 2^{-N}M_d^N$ are
positive for sufficiently large $N$. Thus $A_d$ is {\it 
  irreducible} (see \cite{minc}, Ch.1), so it has a simple
eigenvalue of 1, and all its  other eigenvalues are inside the unit
disk. It follows that $M_d$ has a simple eigenvalue of 2. Let
\begin{equation}
f_d(T) = T^k + c_{k-1}T^{k-1} + \cdots + c_0
\end{equation}
be the minimal polynomial of $M_d$. Let $\rho_d < 2$
be the maximum modulus of any non-2 root of  $f_d$, and let $1+\sigma_d$ be
the maximum multiplicity of any such maximal root. 
 Then for $r \ge 0$ and all $(\al,\be)$,
\begin{equation}
m_{\al\be}^{r+k} +  c_{k-1}m_{\al\be}^{r+k-1} + \cdots +
c_0m_{\al\be}^{r} = 0.
\end{equation}

It follows from the standard theory of linear recurrences that for
some constants $c_{\al\be}$, 
\begin{equation}
m_{\al\be}^{r} = c_{\al\be} 2^r +  \mathcal(r^{\sigma_d}\rho_d^r)\qquad
\text{as } r\to \infty. 
\end{equation}
In particular, $\lim_{r\to\infty} A_d^r = A_{d0}:= [c_{\al\be}]$, and
since $A_d^{r+1} 
= A_d A_d^r$, it follows that each column of $A_{d0}$ is an eigenvector
of $A_d$, corresponding to $\la = 1$. Such eigenvectors are constant
vectors and since $A_{d0}$ is doubly stochastic, we may conclude that
for all $(\al,\be)$, $c_{\al\be} = \frac 1{N_d}$. 
Then there exists $c_d > 0$ so that for $r \ge 0$ and all $(\al,\be)$,
\begin{equation}
\left\vert m_{\al\be}^{r} - \frac{2^r}{N_d}\right\vert <
c_dr_d^{\sigma_d}\rho_d^r. 
\end{equation}
Computations show that for 
for small values of $d$ at least,  $\rho_d = \frac 12$
and $\sigma_d = 0$. In any event, by choosing $2> \bar\rho_d > \rho_d$ if
$\sigma_d > 0$, we can replace $r_d^{\sigma_d}\rho_d^r$ by
$\bar\rho_d^r$ in the upper bound.
Putting this together, we have proved the following theorem.
\begin{theorem}
There exist constants $c_d$ and $\bar\rho_d < 2$ so that if $m \in \mathbb N$
and $\al \in \mathcal S_d$, then for all $r \ge 0$,
\begin{equation}
\left\vert B(\al;2^rm,2^r(m+1)) - \frac {2^r}{N_d} \right\vert < c_d
\bar\rho_d^r.
\end{equation}
\end{theorem}

We now use this result on blocks of length $2^r$ to get our main theorem.
\begin{theorem}
For fixed $d \ge 2$, there exists $\tau_d< 1$ so that, for all $\al \in
\mathcal S_d$, 
\begin{equation}
B(\al;0,N) = \frac N {N_d} + \mathcal O(N^{\tau_d}).
\end{equation}
\end{theorem}
\begin{proof}
By (2.25), we have
\begin{equation}
B(\al;0,N) = \sum_{j=0}^{v-1} B(\al;N_j,N_{j+1}) = 
\sum_{j=1}^{v} B(\al;2^{r_j}M_j,2^{r_j}(M_j+1)).
\end{equation}
It follows that
\begin{equation}
\left\vert B(\al;0,N) - \frac N {N_d} \right \vert \le
c_d(\bar\rho_d^{r_1} 
+ \cdots + \bar\rho_d^{r_v}).
\end{equation}
If $\bar\rho_d \le 1$, the upper bound is $\mathcal O(r_1) = \mathcal O(\log
N) = \mathcal O(N^{\epsilon})$ for any $\ep > 0$.
 If $1 \le \bar\rho_d < 2$, the upper bound is $\mathcal 
O(\bar\rho_d^{r_1}) = \mathcal O(N^{\tau_d})$ for $\tau_d = \frac{\log
  \bar\rho_d}{\log 2}$, since $N \le 2^{r_1+1}$.
\end{proof}
Using the notation (1.11), we have
\begin{equation}
T(N;d,i) = \sum_{\al = (i,j) \in \mathcal S_d} B(\al;0,N),
\end{equation}
and the following is an immediate consequence of Lemma 4.2 and Theorem 4.7.
\begin{corollary}
Suppose $d \ge 2$. Then
\begin{equation}
T(N;d,i) = r_{d,i}N + \mathcal O(N^{\tau_d}),
\end{equation}
where, recalling that $d = \prod p_\ell^{e_\ell}$,
\begin{equation}
r_{d,i} = \frac 1d
\cdot \prod_{p_\ell | i} \frac {p_\ell}{p_\ell+1}  
\cdot \prod_{p_\ell \nmid i} \frac {p_\ell^2}{p_\ell^2-1}.
\end{equation}
\end{corollary}
For example, if $p$ is prime, then $f(p,0) = \frac 1{p+1}$ and $f(p,i)
= \frac p{p^2-1}$ when $p \nmid i$. 

In some sense, the model here is a Markov Chain, if we imagine
going from $m$ to 
$2m$ or $2m+1$ with equal probability, so that the
$B(\be;2^rm,2^r(m+1))$'s represent the distribution of destinations
after $r$ steps. Ken Stolarsky has pointed out that \cite{sch} is
a somewhat different
application of the limiting theory of Markov Chains in a number
theoretic setting.

\section{Small values of $d$}

 It is immediate to see (and to prove) that $2 \ | \ s(n)$ if and only
 if $3 \ | \ 
n$, thus $S_2(n)$ cycles among $\{(0,1), (1,1), (1,0)\}$ and $\tau_2 = 0$.
This generalizes to a family of partition
sequences. Suppose $d \ge 2$ is fixed,
and let 
$b(d;n)$ denote the number of ways that $n$ can be written in the form
\begin{equation}
n = \sum_{i\ge 0} \epsilon_i 2^i,\qquad \epsilon_i \in \{0,\dots,d-1\},
\end{equation}
so that $b(2;n) = 1$. It is shown in \cite{rez1} that 
\begin{equation}
\sum_{n=0}^\infty s(n)X^n = X\prod_{j=0}^\infty\left(1
  + X^{2^j} + X^{2^{j+1}}\right).
\end{equation}
A standard partition argument shows that 
\begin{equation}
\sum_{n=0}^\infty b(d;n)X^n = \prod_{j=0}^\infty \frac {1 - X^{d\cdot
    2^j}}{1 - X^{2^j}}.
\end{equation}
Thus,  $s(n) = b(3;n-1)$. An examination of the product in
(5.3) modulo 2 shows that $b(d;n)$ is odd if and
only if $n \equiv 0,1 \text{ mod } d$ (see  \cite{rez1},
Theorems 5.2 and 2.14.) 

Suppose now that $d=3$.
Write the 8 elements of $\mathcal S_3$ in lexicographic order:
\begin{equation}
(0,1), (0,2), (1,0), (1,1), (1,2), (2,0), (2,1), (2,2).
\end{equation}
Then in the notation of the last section,
\begin{equation}
M_3 = \begin{pmatrix}
1&0&0&1&0&0&0&0 \\ 0&1&0&0&0&0&0&1 \\0&0&1&1&0&0&0&0 \\
0&0&0&0&1&0&1&0 \\ 0&1&1&0&0&0&0&0 \\ 0&0&0&0&0&1&0&1 \\
1&0&0&0&0&1&0&0 \\ 0&0&0&0&1&0&1&0
\end{pmatrix}\ .
\end{equation}
The minimal polynomial of $M_3$ is
\begin{equation}
f_3(T) =T^5 -2T^4+T^3-4T^2+4T = T(T-1)(T-2)(T - \mu)(T- \bar\mu),
\end{equation}
where
\begin{equation}
 \mu = \frac {-1 + \sqrt 7 i}2, \qquad
\bar \mu = \frac {-1 -\sqrt 7 i}2.
\end{equation}
Since the roots of $f_3$ are distinct, we see that for each
$(\al,\be) \in 
\mathcal S_3$, for $r \ge 1$, there exist constants $v_{\al\be i}$ so that
\begin{equation}
m_{\al\be}^{(r)} = v_{\al\be 1} +   v_{\al\be 2}\mu^r +
v_{\al\be 3}\bar\mu^r + \frac 18\cdot 2^r = \frac 18\cdot 2^r + \mathcal
O(2^{r/2}).
\end{equation}
(As it happens, there are only eight distinct sequences 
$m_{\al\be}^{(r)}$.) Corollary 4.8 then implies that
\begin{equation}
\begin{gathered}
T(N;3,0) = \frac N4 + \mathcal O(\sqrt N), \\
T(N;3,1) = \frac {3N}8 + \mathcal O(\sqrt N), \ 
T(N;3,2) = \frac {3N}8 + \mathcal O(\sqrt N).
\end{gathered}
\end{equation}
Since $T(N;3,0) + T(N;3,1) + T(N;3,2) = N$, we gain complete
information from studying $T(N;3,0)$ and
\begin{equation}
\Delta(N) = \Delta_3(N) :=  T(N;3,1) - T(N;3,2).
\end{equation}
(That is, $\Delta_3(N+1) - \Delta_3(N)$ equals $0, 1, -1$ when $s(N)
\equiv 0, 1, 2$ mod 3, respectively.)

To study $T(N;3,0)$, we first define
the set $A_3 \subset \mathbb N$ recursively by:
\begin{equation}
0,5,7 \in A_3, \qquad 0 < n \in A_3 \implies 2n, 8n\pm 5, 8n \pm 7 \in A_3.
\end{equation}
Thus,
\begin{equation}
A_3 = \{ 0, 5, 7, 10,
14, 20, 28, 33, 35, 40, 45, 47, 49, 51, 56, 61, 63,\dots \}.
\end{equation}
 
\begin{theorem}
If $n \ge 0$, then $3 \ |\ s(n)$ if and only if $n \in A_3$.
\end{theorem}
\begin{proof}
It follows recursively from (1.2) or directly from (2.7) that
\begin{equation}
s(2n) = s(n),\quad s(8n\pm 5) = 2s(n) + 3s(n\pm 1), \qquad  s(8n\pm 7)
= s(n) + 3s(n\pm 1). 
\end{equation}
Thus, 3 divides $s(n)$ if and only if 3 divides $s(2n), s(8n\pm
5)$ or $ s(8n\pm 7)$. Since every $n > 1$ can be written uniquely as
$2n', 8n'\pm 5$ or $8n'\pm 7$ with $0 \le n' < n$, the description of 
$A_3$ is complete. 
\end{proof}

In the late 1970's, E. W. Dijkstra \cite{dijk}(pp. 215--6, 230--232) 
studied the Stern sequence under the name ``fusc'', and gave a
different description of $A_3$ (p. 232):
\smallskip
\begin{quote}
Inspired by a recent exercise of Don Knuth I tried to characterize the
arguments $n$ such that $3\ | \ \it{fusc}(n)$. With braces used to denote
zero or more instances of the enclosed, the vertical bar as the BNF
`or', and the question mark `?' to denote either a 0 or a 1, the
syntactical representation for such an argument (in binary) is 
\{0\}1\{?0\{1\}0$\vert$?1\{0\}1\}?1\{0\}. I derived this by considering
-- as a direct derivation of my program -- the finite state automaton
that computes \it{fusc} $(N)$ mod 3.
\end{quote} 
\smallskip
Let
\begin{equation}
a_r = | \{ n \in A_3: 2^r \le n < 2^{r+1} \} | = T(2^{r+1};3,0) -
T(2^r;3,0). 
\end{equation}
It follows from (5.12) that
\begin{equation}
a_0 = a_1=0, \quad a_2=a_3=a_4 = 2, \quad a_5 = 10.
\end{equation}
\begin{lemma}
For $r \ge 3$,  $(a_r)$ satisfies the recurrence 
\begin{equation}
a_r = a_{r-1} + 4a_{r-3}.
\end{equation}
\end{lemma}
\begin{proof}
This is evidently true for $r = 3, 4,5$. If $2^r \le n
< 2^{r+1}$ and $n=2n'$, then $2^{r-1} \le n' < 2^r$, so the even
elements of $A_3$ counted in $a_r$ come from elements of $A_3$ counted in
$a_{r-1}$. If $2^r \le n
< 2^{r+1}$ and $n=8n'\pm5$ or $n=8n'\pm7$, then 
$2^{r-3} < n' < 2^{r-2}$ and $n' \in A_3$. Thus the odd
elements of $A_3$ counted in $a_r$ come (in fours) from elements of
$A_3$ counted in $a_{r-3}$. 
\end{proof}
The characteristic polynomial of the recurrence (5.16) is $T^3-T^2-4$
(necessarily a factor of $f_3(T)$), and has roots $T=2,\ \mu$
and $\bar\mu$. The details of the following routine computation are omitted.
\begin{theorem}
For $r \ge 0$, we have the exact formula
\begin{equation}
a_r = \frac 14 \cdot 2^r +  \left(\frac{-7+5\sqrt 7 i}{56}\right)\mu^r +
 \left(\frac{-7-5\sqrt 7 i}{56}\right)\bar\mu^r .
\end{equation}
\end{theorem}
Keeping in mind that $s(0) = 0$ is not counted in any $a_r$, we find
 after a further computation that the error estimate $\mathcal O(\sqrt
 N)$ is best possible for $T(N;3,0)$: 
\begin{corollary}
\begin{equation}
T(2^r;3,0) = \frac 14 \cdot 2^r +  \left(\frac{7-\sqrt 7
  i}{56}\right)\mu^r +  \left(\frac{7+\sqrt 7 i}{56}\right)\bar\mu^r
+ \frac 12\ .
\end{equation}
\end{corollary}

To study $\Delta(N)$, we first need a somewhat surprising lemma.
\begin{lemma}
For all $N$, $\Delta(2N)= \Delta(4N)$. 
\end{lemma}
\begin{proof}
The simplest proof is by induction, and the assertion is trivial for $N=0$.
 There are eight possible ``short'' diatomic arrays modulo 3:
\begin{equation}
\begin{gathered}
\begin{matrix}
  {s(N)}&&&&{s(N+1)} \\  s(2N)&&s(2N+1)&&s(2N+2) \\
  s(4N)&s(4N+1)&s(4N+2)&s(4N+3)&s(4N+4) 
\end{matrix} = \\
\begin{matrix}
0&&&&1\\0&&1&&1\\0&1&1&2&1
\end{matrix}\quad\Bigg\vert\Bigg\vert\quad 
\begin{matrix}
0&&&&2\\0&&2&&2\\0&2&2&1&2
\end{matrix}\quad \Bigg\vert\Bigg\vert\quad
\begin{matrix}
1&&&&0\\1&&1&&0\\1&2&1&1&0
\end{matrix}\quad\Bigg\vert\Bigg\vert\quad 
\begin{matrix}
1&&&&1\\1&&2&&1\\1&0&2&0&1
\end{matrix}  \\
\begin{matrix}
1&&&&2\\1&&0&&2\\1&1&0&2&2
\end{matrix} \quad\Bigg\vert\Bigg\vert\quad 
\begin{matrix}
2&&&&0\\2&&2&&0\\2&1&2&2&0
\end{matrix} \quad\Bigg\vert\Bigg\vert\quad 
\begin{matrix}
2&&&&1\\2&&0&&1\\2&2&0&1&1
\end{matrix} \quad\Bigg\vert\Bigg\vert\quad 
\begin{matrix}
2&&&&2\\2&&1&&2\\2&0&1&0&2
\end{matrix}
 \end{gathered}
\end{equation}
By counting the elements in the rows mod 3 in each case, 
we see that  $\Delta(2N+2) - \Delta(2N) =
\Delta(4N+4) - \Delta(4N)$ is equal to:  $1, -1, 2, 0, 1,
-2, -1, 0$, respectively.
\end{proof}
\begin{theorem}
For all $n$, $\Delta(n) \in \{0,1,2,3\}$. More specifically,
\begin{equation}
\begin{gathered}
S_3(m) = (0,1) \implies \Delta(2m) = 0,\ \Delta(2m+1) =0; \\ 
S_3(m) = (0,2) \implies \Delta(2m) = 3,\ \Delta(2m+1) =3; \\ 
S_3(m) = (1,*) \implies \Delta(2m) = 1,\ \Delta(2m+1) =2; \\ 
S_3(m) = (2,*) \implies \Delta(2m) = 2,\ \Delta(2m+1) =1. 
\end{gathered}
\end{equation}
\end{theorem}
\begin{proof}

To prove the theorem, we first observe that it is correct for $m \le
4$. We now assume it is true for $m \le
2^r$ and prove it for $2^r \le m < 2^{r+1}$. There are sixteen
cases: $m$ can be even or odd and there are 8 choices for $S_3(m)$.
As a representative example, suppose $S_3(m) = (2,1)$. We shall
consider the cases $m=2t$ and $m=2t+1$ separately. The proofs for the
other seven choices of $S_3(m)$ are very similar and are omitted.

Suppose first that $m=2t < 2^{r+1}$. Then $S_3(m) = S_3(2t) = (2,1)$, hence
$S_3(t) = (2,2)$. We have $\Delta(2m) = 2$ by hypothesis, and hence
$\Delta(4m) = 2$ by Lemma 5.5. The eighth array in (5.19) shows that
$s(4t) \equiv 2$ mod 3, so that $\Delta(4m+1) = \Delta(4m) - 1 = 1$,
as asserted in (5.20).

If, on the other hand, $m=2t+1 < 2^{r+1}$ and  $S_3(m) = S_3(2t+1) = (2,1)$,
then $S_3(t) = (1,1)$. We now have $\Delta(2t) = 1$ and $\Delta(2t+1)
=2$ by hypothesis and $\Delta(4t) = 1$ by Lemma 5.5. The fourth array
in (5.19) shows that $(s(4t),s(4t+1),s(4t+2)) \equiv (1,0,2)\text{ mod 
}3$. Thus, it follows that $\Delta(2m) = \Delta(4t+2) = \Delta(4t) + 1
+ 0 =2$ and 
 $\Delta(2m+1) = \Delta(4t+3) =  \Delta(4t+2)-1 = 1$, again as desired.
\end{proof}

Since $S_3(m)$ is uniformly distributed on $\mathcal S_3$, (5.20)
shows that $\Delta(n)$ takes the values $(0,1,2,3)$ with limiting
probability $(\frac 18, \frac 38, \frac 38, \frac 18)$.

We conclude with a few words about the results announced at the end of the
first section, but not proved here. For each $(d,i)$, $T(2^r;d,i)$ will
satisfy a recurrence whose characteristic equation is a factor of the
minimal polynomial of $\mathcal S_d$. It happens that $T(2^r;4,0) =
T(2^r;5,0)$ for small values of $r$ and both 
satisfy the recurrence with characteristic polynomial
$T^4-2T^3+T^2-4$ (roots: $2,-1,-\tau, -\bar\tau$)
so that equality holds for all $r$. The same
applies to  $T(2^r;6,0) = T(2^r;9,0) = T(2^r;11,0)$, with a more
complicated recurrence.
Results similar to Lemma 5.5 and Theorem 5.6 hold for $d=4$, with a
similar proof; 
Antonios Hondroulis has shown that this is also true for $d=6$. No result
has been found yet for $d=5$, although a Mathematica check for $N \le
2^{19}$ shows that $-5 \le T(N;5,1) - T(N;5,4) \le 11$. These topics
will be discussed in greater detail in \cite{rez2}.

\bibliographystyle{amsalpha}

\begin{thebibliography}{99}

\bibitem{as} J.-P. Allouche, J. Shallit, The ring of $k$-regular
  sequences, {\it Theoret. Comput. Sci.} {\bf 98} (1992), 163--197,
  MR1166363 (94c:11021).

\bibitem{bach} P. Bachmann, Niedere Zahlentheorie, v. 1, Leipzig 1902,
  Reprinted by Chelsea, New York, 1968.

\bibitem{broc} A. Brocot, Calcul des rouages par approximation,
  nouvelle m\'ethode, {\it Revue chronom\'etrique. Journal des
    horlogers, scientifique et pratique} {\bf 3} (1861), 186--194.

\bibitem{cw} N. Calkin and H. Wilf, Recounting the rationals, {\it
  Amer. Math. Monthly} {\bf 107} (2000), 360--363, MR1763062 (2001d:11024).

\bibitem{dr} G. de Rham,  Un peu de math\'ematiques \`a propos d'une
  courbe plane. {\it Elemente der Math.} {\bf 2} (1947), 73--76, 89--97
  MR0022685 (9,246g), reprinted in Oeuvres Math\'ematiques, Geneva,
  1981, 678--690, MR0638722 (84d:01081).

\bibitem{dic} L. E. Dickson,  History of the Theory of Numbers, v. 1,
 Carnegie Inst. of Washington, Washington, D.C., 1919, reprinted by
 Chelsea, New York, 1966,  MR0245499 (39 \#6807a).

\bibitem{dijk} E. W. Dijkstra,  Selected writings on computing: a
  personal perspective, Springer-Verlag, New York, 1982,
MR0677672 (85d:68001).

\bibitem{gkp} R. L. Graham, D. E. Knuth, O. Patashnik, Concrete
  Mathematics, Second Edition, Addison-Wesley, Boston, 1994, MR1397498
  (97d:68003). 

\bibitem{hayes} B. Hayes, On the teeth of wheels, {\it American
  Scientist} {\bf 88},  July-August 2000, 296--300.

\bibitem{hkmpp} A. Hinz, S. Klav\v zar, U. Milutinovi\'c, D. Parisse,
C. Petr, Metric properties of the Tower of Hanoi graphs and Stern's
diatomic sequence, {\it Euro. J. Comb.} {\bf 26} (2005), 693--708,
MR2127690 (2005m:05081).

\bibitem{leh} D. H. Lehmer, On Stern's diatomic series, {\it
  Amer. Math. Monthly}  {\bf 36} (1929), 59--67, MR1521653.

\bibitem{luc} E. Lucas, Th\'eorie des nombres, vol. 1, Gauthier-Villars,
  Paris, 1891.

\bibitem{minc} H. Minc, Nonnegative matrices, Wiley, New York, 1988,
  MR0932967 (89i:15061).

\bibitem{mink} H. Minkowski, Zur Geometrie der Zahlen, Ver. III
  Int. Math.-Kong. Heidelberg 1904, pp. 164-173; in {\it
  Gesammelte Abhandlungen, Vol. 2}, Chelsea, New York
  1967, pp. 45--52.

\bibitem{rez1} B. Reznick, Some digital partition functions, in:
  B.C. Berndt et al. (Eds.), Analytic Number Theory, Proceedings of a
  Conference in Honor of Paul T. Bateman, Birkh\"auser, Boston, 1990,
  pp. 451--477, MR1084197 (91k:11092).

\bibitem{rez2} B. Reznick, A Stern introduction to combinatorial number
  theory, in preparation.

\bibitem{sch} W. Schmidt, The joint distribution of the digits of
  certain integer $s$-tuples, in: Paul Erd\"os, Editor-in-Chief,
  Studies in Pure Mathematics to the memory of Paul Tur\'an,
  Birkh\"auser, Basel, 1983, pp. 605--622, MR0820255 (87h:11072).

\bibitem{ste} M. A. Stern, Ueber eine zahlentheoretische Funktion,
  {\it J. Reine Angew. Math.} {\bf 55} (1858) 193--220.

\bibitem{urb} I. Urbiha,  Some properties of a function studied by de
  Rham, Carlitz and Dijkstra and its relation to the
  (Eisenstein-)Stern's diatomic sequence, {\it Math. Commun} {\bf 6}
  (2001) 181--198 (MR1908338 (2003f:11018) 

\end{thebibliography}

\end{document}